\documentclass[a4paper,12pt]{amsart}
\usepackage{amssymb,amsmath,a4wide,graphicx,stmaryrd,fullpage,setspace,microtype,xr}
\usepackage[pagebackref=true]{hyperref}
\usepackage[all]{xy}
\title{Flows revisited: the model category structure and its left determinedness}

\author[P. Gaucher]{Philippe Gaucher}

\address{Universit\'e de Paris, IRIF, CNRS, F-75013 Paris, France}

\urladdr{http://www.irif.fr/{\~{}}gaucher} 



\swapnumbers

\newcommand{\C}{\mathcal{C}}

\newcommand{\K}{\mathcal{K}}

\newcommand{\W}{\mathcal{W}}

\newcommand{\p}\times

\renewcommand{\P}{\mathbb{P}}

\newtheorem{thm}{Theorem}[section]
\newtheorem{prop}[thm]{Proposition}
\newtheorem{lem}[thm]{Lemma}

\newtheorem{defn}[thm]{Definition}
\newtheorem{nota}[thm]{Notation}

\newcommand{\bd}{\begin{defn}}
\newcommand{\ed}{\end{defn}}
\newcommand{\bp}{\begin{prop}}
\newcommand{\ep}{\end{prop}}
\newcommand{\bth}{\begin{thm}}
\renewcommand{\eth}{\end{thm}}
\newcommand{\bpf}{\begin{proof}}
\newcommand{\epf}{\end{proof}}

\newcommand{\fL}[1]{\ar@{->}[ll]_-{#1}}
\newcommand{\fR}[1]{\ar@{->}[rr]^-{#1}}
\newcommand{\fRr}[1]{\ar@{->}[rrr]^-{#1}}
\newcommand{\fD}[1]{\ar@{->}[dd]_-{#1}}
\newcommand{\fU}[1]{\ar@{->}[uu]^-{#1}}
\newcommand{\f}[2]{\ar@{->}[#1]|{#2}}
\newcommand{\ff}[2]{\ar@2{->}[#1]|{#2}}
\newcommand{\frr}[1]{\ar@{->}[rrrr]^-{#1}}

\newcommand{\fl}[1]{\ar@{->}[l]_-{#1}}
\newcommand{\fr}[1]{\ar@{->}[r]^-{#1}}
\newcommand{\fd}[1]{\ar@{->}[d]_-{#1}}
\newcommand{\fu}[1]{\ar@{->}[u]^-{#1}}

\renewcommand{\top}{{\mathbf{Top}}}

\newcommand{\iso}{\cong}

\renewcommand{\leq}{\leqslant}
\renewcommand{\geq}{\geqslant}

\def\cartesien{%
  \ar@{-}[]+R+<6pt,-2pt>;[]+RD+<6pt,-6pt>%
  \ar@{-}[]+D+<2pt,-6pt>;[]+RD+<6pt,-6pt>%
}
\def\cocartesien{%
  \ar@{-}[]+L+<-6pt,+2pt>;[]+LU+<-6pt,+6pt>%
  \ar@{-}[]+U+<-2pt,+6pt>;[]+LU+<-6pt,+6pt>%
}
\def\hocartesien{%
  \ar@{-}[]+R+<6pt,-2pt>;[]+RD+<6pt,-6pt>_{h}%
  \ar@{-}[]+D+<2pt,-6pt>;[]+RD+<6pt,-6pt>%
}
\def\hococartesien{%
  \ar@{-}[]+L+<-6pt,+2pt>;[]+LU+<-6pt,+6pt>_{h}%
  \ar@{-}[]+U+<-2pt,+6pt>;[]+LU+<-6pt,+6pt>%
}

\newcommand{\brm}[1]{\rm{\mathbf{#1}}}

\newcommand{\dtop}{{\brm{Flow}}}
\newcommand{\set}{{\brm{Set}}}

\newcommand{\ttop}{{\brm{TOP}}}

\newcommand{\glob}{{\rm{Glob}}}

\DeclareMathOperator{\id}{Id}

\DeclareMathOperator{\cocyl}{{Path}}

\newcommand{\liminj}{\varinjlim}
\newcommand{\limproj}{\varprojlim}

\makeatletter
\def\varholim@#1#2{%
  \vtop{\m@th\ialign{##\cr
    \hfil$#1\operator@font holim$\hfil\cr
    \noalign{\nointerlineskip\kern1.5\ex@}#2\cr
    \noalign{\nointerlineskip\kern-\ex@}\cr}}%
}
\def\holimproj{%
  \mathop{\mathpalette\varholim@{\leftarrowfill@\textstyle}}\nmlimits@
}
\def\holiminj{%
  \mathop{\mathpalette\varholim@{\rightarrowfill@\textstyle}}\nmlimits@
}
\makeatother

\DeclareMathOperator{\cell}{{\brm{cell}}}
\DeclareMathOperator{\cof}{{\brm{cof}}}
\DeclareMathOperator{\inj}{{\brm{inj}}}

\begin{document}

\begin{abstract} 
	Flows are a topological model of concurrency which enables to encode the notion of refinement of observation and to understand the homological properties of branchings and mergings of execution paths. Roughly speaking, they are Grandis' $d$-spaces without an underlying topological space. They just have an underlying homotopy type. This note is twofold. First, we give a new construction of the model category structure of flows which is more conceptual thanks to Isaev's results. It avoids the use of difficult topological arguments. Secondly, we prove that this model category is left determined by adapting an argument due to Olschok. The introduction contains some speculations about what we expect to find out by localizing this minimal model category structure.
	\linebreak[4]\linebreak[4]
	Les flots sont un mod\`ele topologique de la concurrence qui permet d'encoder la notion de raffinement de l'observation et de comprendre les propri\'et\'es homologiques des branchements et des confluences des chemins d'ex\'ecution. Intuitivement, ce sont des $d$-espaces au sens de Grandis sans espace topologique sous-jacent. Ils ont seulement un type d'homotopie sous-jacent. Cette note a deux objectifs. Premi\`erement de donner une nouvelle construction de la cat\'egorie de mod\`eles des flots plus conceptuelle gr\^ace au travail d'Isaev. Cela permet d'\'eviter des arguments topologiques difficiles. Deuxi\`emement nous prouvons que cette cat\'egorie de mod\`eles est d\'etermin\'ee à gauche en adaptant un argument de Olschok. L'introduction contient quelques sp\'eculations sur ce qu'on s'attend \`a trouver en localisant cette cat\'egorie de mod\`eles minimale.
\end{abstract}

\maketitle
\tableofcontents

\section{Introduction}

\subsection{Topological models of concurrency}
There is a multitude of topological models of concurrency: flows which are the subject of this paper and which are introduced in \cite{model3}, but also $d$-spaces \cite{mg}, streams \cite{MR2545830}, inequilogical spaces \cite{equilogical}, spaces with distinguished cu\-bes \cite{distinguishedcube1} \cite{distinguishedcube2}, multipointed $d$-spaces \cite{mdtop} etc... All these mathematical devices contain the same basic examples coming from concurrency (e.g. the geometric realizations of precubical sets), a local ordering modeling the direction of time and its irreversibility, execution paths, a set or a topological space of states and a notion of homotopy between execution paths to model concurrency. Grandis' $d$-spaces give rise to a vast literature studying the directed fundamental category and the directed components of directed spaces whatever the definition we give to this notion of \textit{directed space}. 

This paper belongs to the sequence of papers \cite{model3} \cite{hocont} \cite{1eme} \cite{2eme} \cite{3eme} \cite{4eme} \cite{nonexistence} \cite{exbranch}. The main feature of the model category of flows is to enable the formalization and the study of the notion of refinement of observation (the cofibrant replacement functor plays a crucial role in the formalization indeed). The model category of flows also enables the study of the homological properties of the branching areas and merging areas of execution paths in concurrent systems, in particular a long exact sequence, and their interaction with the refinement of observation, actually their invariance with respect to them. Since flows have also a labeled version (see \cite[Section~6]{ccsprecub}), they can be used for modeling the path spaces of process algebras for any synchronization algebra \cite{ccsprecub} \cite{symcub}. 

\subsection{Speculative digression}

In our line of research, the objects are multipointed, i.e. they are equipped with a distinguished set of states. The set of states provided by the description of a concurrent process is not forgotten. It is exactly the same phenomenon as in the formalism of simplicial set. A simplicial set is equipped with a set of vertices which comes from the description of the space. This family of topological models of concurrency differs from other topological models of concurrency like Grandis' $d$-spaces or streams which are not multipointed. 

The interest of our approach is that it is already possible to build model category structures such that the weak equivalences preserve the causal stru\-cture of a process. Indeed, all flows are fibrant. By Ken Brown's lemma, two flows are therefore weakly equivalent if and only if they can be related by a zig-zag of trivial fibrations. Every trivial fibration satisfies the right lifting property with respect to any cofibration, and therefore with respect to any reasonable notion of extension of paths. Thus two weakly equivalent flows are bisimilar in Joyal-Nielsen-Winskel's sense \cite{0856.68067} for any reasonable notion of extension of paths.

\begin{figure}
	\begin{center}
		\includegraphics[width=8cm]{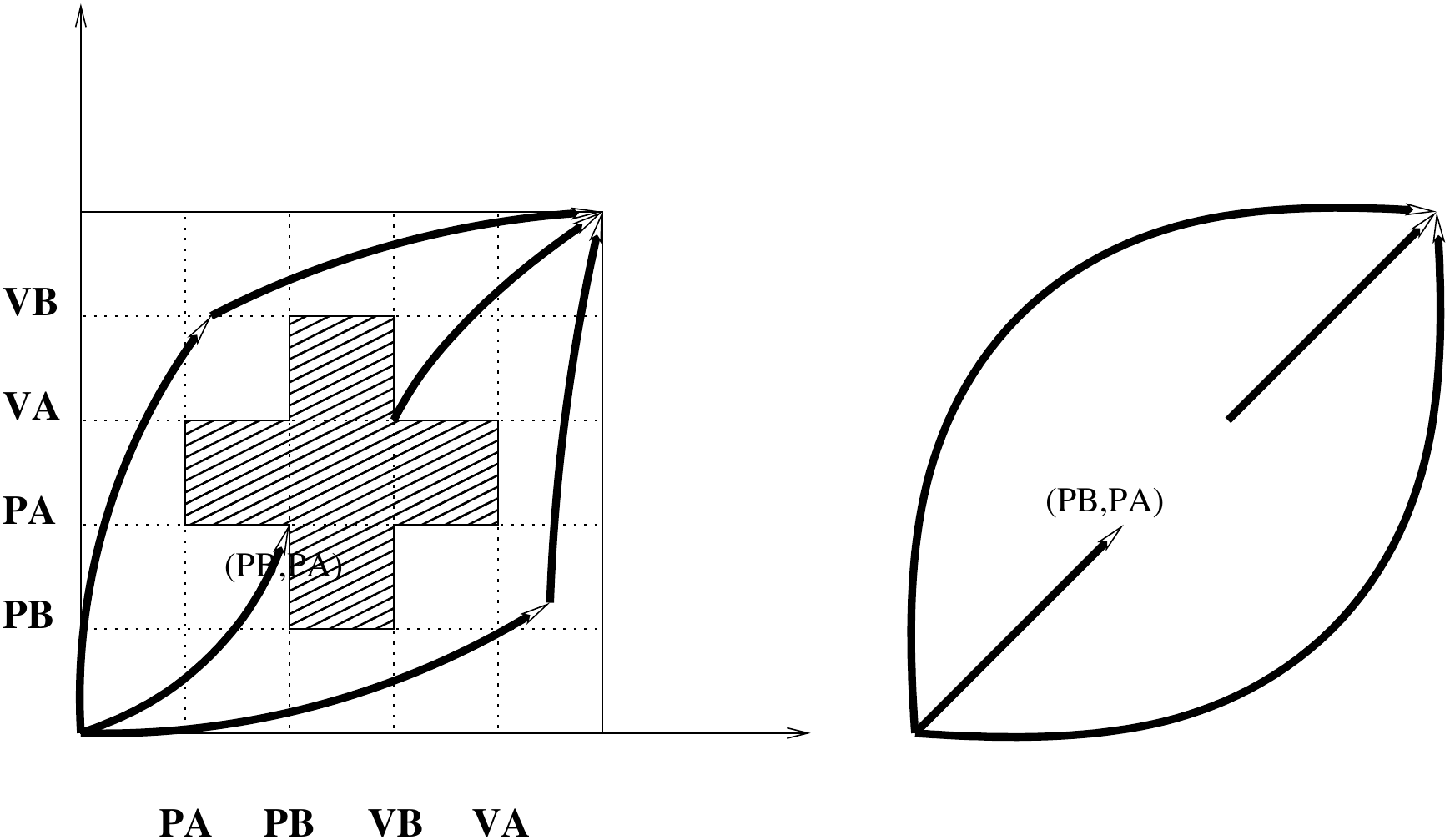}
	\end{center}
	\caption{Swiss Flag example}
	\label{ex3}
\end{figure}

The main drawback of a lot of, and to the best of our knowledge, actually all other model categories introduced in directed homotopy is that their weak equivalences destroy the causal structure for example by identifying the directed segment up to weak equivalence with a point. Indeed, the directed segment should not be contractible (in a directed sense) in directed homotopy. To understand the reason, consider the well-known example of the Swiss Flag example (cf. Figure~\ref{ex3}). It consists of two processes concurrently executing the instructions $PA.PB.VA.VB$ and $PB.PA.VB.VA$ where $Px$ means taking the control of a shared ressource $x$ and $Vx$ means releasing it. We suppose that at most one process can take the control of a given shared ressource. There are in our example two shared ressources $A$ and $B$ which can be for example buffers to temporarily store data. The execution paths start from the bottom left corner, ends to the top right corner, and are supposed to be nondecreasing with respect to each axis of coordinates. Then the point of coordinates $(PB,PA)$ is a deadlock because it is impossible from it to reach the desired final state. As soon as the directed segment becomes weakly equivalent to a point, the deadlock $(PB,PA)$ could disappear up to weak equivalence (the picture on the right is an object having the same causal structure). It means that a weak equivalence could break the causal structure as soon as the directed segment is contractible in a directed sense and, in this case, if left properness is assumed.

\begin{figure}
	\begin{center}
		\includegraphics[width=8cm]{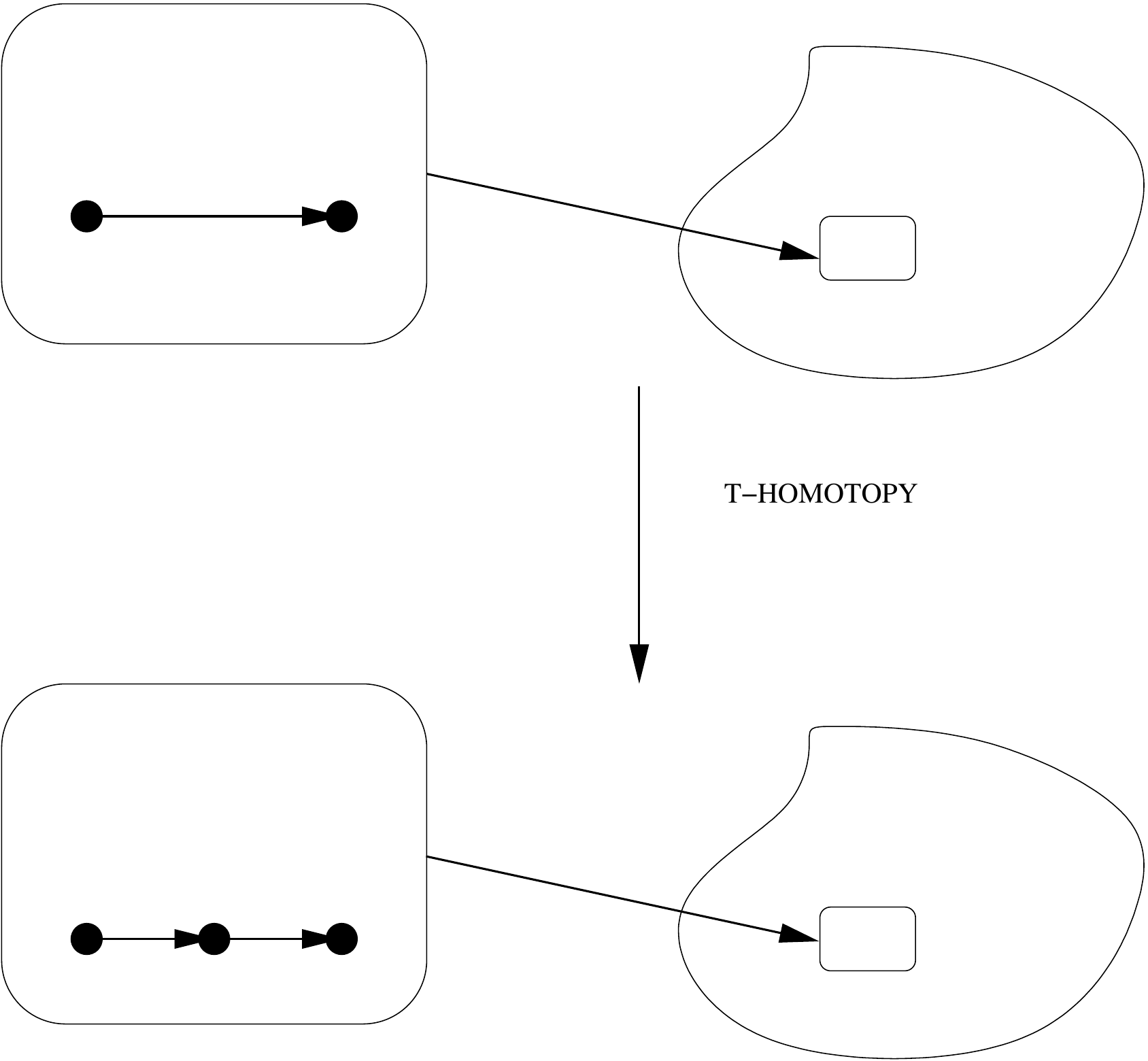}
	\end{center}
	\caption{Replacement of a directed segment by a more refined one}
	\label{ex2}
\end{figure}

The weak equivalences of the model category of flows do preserve the causal structure. However, they are too restrictive. It is not even possible up to weak equivalence to replace in a flow a directed segment by a more refined one (cf. Figure~\ref{ex2}), by adding additional points in the middle of the segment in the distinguished set of states. This annoying behaviour can be overcome by adding weak equivalences by homotopical localization. One of the challenges of our line of research is precisely to understand the homotopical localization of the model category of (labeled) flows with respect to this kind of maps: they are called T-homotopy equivalences in \cite{3eme}. 

The paper \cite{nonexistence} proves not only that the class of weak equivalences of this homotopical localization contains more equivalences than the dihomotopy equivalences of flows as defined in \cite{hocont}, but also that there is no hope to obtain a model category structure on flows such that the weak equivalences are exactly these dihomotopy equivalences, even if a notion of fibrant object (the homotopy continuous flows \cite[Definition~4.3]{hocont}) with the associated Whitehead theorem exists for dihomotopy equivalences \cite[Theorem~4.6]{hocont}. 

By now, we only know by studying examples that the weak equivalences of this homotopical localization seem to be, in the labeled case, dihomotopy equivalences in the sense of \cite{hocont} up to a kind of bisimulation. In particular, it means that the weak equivalences of this homotopical localization likely preserve causality, but not the branching and merging homologies, and not the underlying homotopy type of a flow as defined in \cite{4eme} which is, roughly speaking, the homotopy type of the space obtained after removing the directed structure of a flow. 

It is actually possible to homotopically localize the model category of labeled flows by the whole class of bisimulations in Joyal-Nielsen-Winskel's sense since this class of maps is accessible. Then another strange phenomenon occurs. We would then have to deal in the localization with weak equivalences breaking the causal structure. The latter phenomenon is explained in \cite[Theorem~12.4]{pastsim} within the combinatorial framework of Catta\-ni-Sassone higher dimensional transition systems but it can be easily adapted and generalized to many other frameworks of directed homotopy, including the one of flows.

\subsection{Purpose of this note}

The construction of the model category of flows as carried out in \cite{model3} is quite long and tricky. It makes use of rather complicated topological lemmas, in particular because colimits of flows are difficult to understand. Indeed, flows are roughly speaking \textit{small categories without identities (precategories ? pseudocategories ?) enriched over topological spaces}~\footnote{There is no known relation between the model category of flows and the model category of topologically enriched small categories of \cite{MR3094501}: the obvious adjunction is not even a Quillen adjunction; moreover the terminology of ``flow'' must not mislead the reader, it has nothing to do with a similar terminology in Morse theory.} Therefore new paths are created as soon as states are identified in a colimit, which may generate complicated modifications of the topology of the path space (e.g. see the proof of \cite[Theorem~15.2]{model3} which is not only complicated but also contains a flaw which will be fixed in a subsequent paper). The first purpose of this note is to drastically simplify this construction using \cite{Isaev}. We then explain in a second part why the model category of flows is left determined in the sense of \cite{rotho} by adapting an argument due to Marc Olschok for the model category of topological spaces. The latter fact is a new result (it was mentioned without proof in \cite[Section~12]{pastsim}). In the (complicated) quest of finding out better model categories with as much weak equivalences as possible preserving the causal structure, this result means that it is not possible in the framework of flows to remove weak equivalences without changing the set of generating cofibrations.

\subsection{Organization}

Section~\ref{section1} recalls what we need to use from Isaev's paper.  Section~\ref{section2} explains the new construction of the model category structure of flows (Theorem~\ref{simplification}). Section~\ref{section3} recalls the notion of left determined model category and proves that the category of flows is left determined (Theorem~\ref{leftdetflow}). Section~\ref{conc} makes some final comments.

\subsection{Notations}

All categories are locally small. The category of sets is denoted by $\set$. The set of maps in a category $\K$ from $X$ to $Y$ is denoted by $\K(X,Y)$.  The initial (final resp.) object, if it exists, is always denoted by $\varnothing$ ($\mathbf{1}$ resp.). The identity of an object $X$ is denoted by $\id_X$.  The composite of two maps $f:A\to B$ and $g:B\to C$ is denoted by $g.f$. A subcategory is always isomorphism-closed (replete). Let $f$ and $g$ be two maps of a category $\K$. Denote by $f\boxslash g$ when $f$ satisfies the \emph{left lifting property} (LLP) with respect to $g$, or equivalently $g$ satisfies the \emph{right lifting property} (RLP) with respect to $f$. Let $\C$ be a class of maps. Let us introduce the notations $\inj(\C) = \{g \in \K, \forall f \in \C, f\boxslash g\}$ and $\cof(\C)=\{f\mid \forall g\in \inj(\C), f\boxslash g\}$. The class of morphisms of $\K$ that are transfinite compositions of pushouts of elements of $\C$ is denoted by $\cell(\C)$. We refer to \cite{MR95j:18001} for locally presentable categories, to \cite{MR2506258} for combinatorial model categories.  We refer to \cite{MR99h:55031} and to \cite{ref_model2} for more general model categories.

\section{Isaev approach for constructing model categories}
\label{section1}

Let $\K$ be a locally presentable category. A combinatorial model category structure is characterized by its set of generating cofibrations and by its class of fibrant objects by \cite[Proposition E.1.10]{quasicat}. Therefore, for a given set of maps $I$, there exists at most one combinatorial model category structure on $\K$ such that the set of generating cofibrations is $I$ and such that all objects are fibrant. In \cite{Isaev}, several methods are expounded to obtain model category structures such that all objects are fibrant. We summarize in the next theorem what we are going to need in this note. 

\bth \cite[Theorem~4.3, Proposition~4.4, Proposition~4.5 and Corollary~4.6]{Isaev} \label{Isa}
Let $\K$ be a locally presentable category. Let $I$ be a set of maps of $\K$ such that the domains of the maps of $I$ are $I$-cofibrant (i.e. belong to $\cof(I)$). Suppose that for every map $i:U\to V \in I$, the relative codiagonal map $V\sqcup_U V \to V$ factors as a composite $V\sqcup_U V \stackrel{\gamma_0\sqcup \gamma_1}\to C_U(V)\to V$ such that the left-hand map belongs to $\cof(I)$. Let $J_I = \{\gamma_0:V\to C_U(V)\mid U\to V \in I\}$. Suppose that there exists a path functor $\cocyl:\K \to \K$, i.e. an endofunctor of $\K$ equipped with two natural transformations $\tau:\id\Rightarrow\cocyl$ and $\pi:\cocyl\Rightarrow \id\p \id$ such that the composite $\pi.\tau$ is the diagonal. Moreover we suppose that the path functor satisfies the following hypotheses: 
\begin{enumerate}
	\item With $\pi=(\pi_0,\pi_1)$, $\pi_0:\cocyl(X)\to X$ and $\pi_1:\cocyl(X)\to X$ have the RLP with respect to $I$.
	\item The map $\pi:\cocyl(X)\to X\p X$ has the RLP with respect to the maps of $J_I$.
\end{enumerate}
Then there exists a unique model category structure on $\K$ such that the set of generating cofibrations is $I$ and such that the set of generating trivial cofibrations is $J_I$. Moreover, all objects are fibrant.
\eth

Unlike in \cite{Isaev}, we can drop the hypothesis about the smallness of the domains and the codomains of the maps of $I$ with respect to $I$ by \cite[Proposition~1.3]{MR1780498} because the ambient category is supposed to be locally presentable . Note that every map of $J_I$ is a split monomorphism since the composite $V\sqcup_U V \to C_U(V)\to V$ is the relative codiagonal. Therefore every object is fibrant indeed.

\section{The model category of flows}
\label{section2}

\begin{nota} The category $\top$ denotes a bicomplete locally presentable cartesian closed full subcategory of the category of general topological spa\-ces containing all CW-complexes. \end{nota}

The category of $\Delta$-generated spaces, i.e. the colimits of simplices, or equivalently the colimits of the segment $[0,1]$ by \cite[Proposition~3.17]{LisDicov}, satisfies these hypotheses \cite{FR}. It is also possible to add weak separability hypotheses like this one: for every continuous map $g:\Delta^n\to X$ where $\Delta^n$ is the topological $n$-simplex with $n\geq 0$, $g(\Delta^n)$ is closed in $X$. For a tutorial about these topological spaces, see for example \cite[Section~2]{mdtop}. We take $\top$ to be equipped with the standard Quillen model category structure. 

\begin{nota}
	The internal hom functor is denoted by $\ttop(-,-)$.
\end{nota}

\bd \cite{model3} A {\rm flow} $X$ consists of a topological space $\P X$ of execution paths, a discrete space $X^0$ of states, two continuous maps $s$ and $t$ from $\P X$ to $X^0$ called the source and target map respectively, and a continuous and associative map \[*:\{(x,y)\in \P X\p \P X; t(x)=s(y)\}\longrightarrow \P X\] such that $s(x*y)=s(x)$ and $t(x*y)=t(y)$.  A morphism of flows $f:X\longrightarrow Y$ consists of a set map $f^0:X^0\longrightarrow Y^0$ together with a continuous map $\P f:\P X\longrightarrow \P Y$ such that $f(s(x))=s(f(x))$, $f(t(x))=t(f(x))$ and $f(x*y)=f(x)*f(y)$. The corresponding category is denoted by $\dtop$. \ed

\begin{nota}
	For a topological space $X$, let $\glob(X)$ be the flow defined by $\glob(X)^0=\{0,1\}$ and $\P \glob(X)=X$ with $s$ and $t$ being the constant functions $s=0$ and $t=1$. The Glob mapping induces a functor from the category $\top$ of topological spaces to the category $\dtop$ of flows.
\end{nota}

We need to recall the two following easy propositions: 

\bp \label{ortho3} \cite[Proposition~13.2]{model3} A morphism of flows
$f:X\longrightarrow Y$ satisfies the RLP with respect to
$\glob(U)\longrightarrow \glob(V)$ if and only if for any
$\alpha,\beta\in X^0$, $\P_{\alpha,\beta}X\longrightarrow
\P_{f(\alpha),f(\beta)}Y$ satisfies the RLP with respect to
$U\longrightarrow V$. \ep

\bp\label{ortho2} \cite[Proposition~16.2]{model3} Let $f$ be a morphism of flows. Then the
following conditions are equivalent:
\begin{enumerate}
	\item $f$ is bijective on states
	\item $f$ satisfies the RLP
	with respect to $R:\{0,1\}\longrightarrow \{0\}$ and
	$C:\varnothing\subset \{0\}$.
\end{enumerate}
\ep

We will also need this new proposition which does not seem to be proved in one of our previous papers about flows:  

\bp \label{connected-colim-glob}
The globe functor $\glob:\top\to \dtop$ preserves connected colimits (i.e. colimits such that the underlying small category is connected).
\ep

Note that the connectedness hypothesis is necessary. Indeed, $V$ and $W$ being two topological spaces, the flow $\glob(V\sqcup W)$ has two states whereas the flow $\glob(V) \sqcup \glob(W)$ has four states.

\bpf Let $V$ be a topological space. Giving a map from the flow $\glob(V)$ to a flow $X$ is equivalent to choosing two states $\alpha$ and $\beta$ of $X$ (the image of the states $0$ and $1$ of $\glob(V)$) and a continuous map from $V$ to $\P_{\alpha,\beta}X$. Thus the following natural bijection of sets holds \begin{equation}
\label{truc}
\dtop(\glob(V),X)\iso \bigsqcup_{(\alpha,\beta)\in X^0\p X^0} \top(V,\P_{\alpha,\beta}X).
\end{equation}
We obtain the sequence of natural bijections ($\liminj V_i$ being a connected colimit of topological spaces)
\begin{align*}
\dtop(\glob(\liminj V_i),X) & \iso \bigsqcup_{(\alpha,\beta)\in X^0\p X^0} \top(\liminj V_i,\P_{\alpha,\beta}X) \\
&\iso \bigsqcup_{(\alpha,\beta)\in X^0\p X^0} \limproj\top(V_i,\P_{\alpha,\beta}X)\\ 
&\iso \limproj \bigsqcup_{(\alpha,\beta)\in X^0\p X^0} \top(V_i,\P_{\alpha,\beta}X) \\
&\iso \limproj \dtop(\glob(V_i),X) \\
&\iso \dtop(\liminj \glob(V_i),X),
\end{align*}
the first and the fourth isomorphisms by (\ref{truc}), the second and the fifth isomorphisms by definition of a (co)limit and the third isomorphism by the connectedness of the limit. 
The proof is complete using the Yoneda lemma.
\epf

\begin{nota} \cite[Notation~7.6]{model3} Let $U$ be a topological space. Let $X$ be a flow. The flow $\{U,X\}_S$ is defined as follows:
	\begin{enumerate}
		\item The set of states of $\{U,X\}_S$ is $X^0$.
		\item For $\alpha,\beta\in X^0$, let $\P_{\alpha,\beta}\{U,X\}_S = \ttop(U,\P_{\alpha,\beta}X)$.
		\item For $\alpha,\beta,\gamma\in X^0$, the composition law
		\[*:\P_{\alpha,\beta}\{U,X\}_S\p
		\P_{\beta,\gamma}\{U,X\}_S\longrightarrow\P_{\alpha,\gamma}\{U,X\}_S\]
		is the composite 
		\[\P_{\alpha,\beta}\{U,X\}_S\p \P_{\beta,\gamma}\{U,X\}_S\iso
		\ttop\left(U,\P_{\alpha,\beta}X\p
		\P_{\beta,\gamma}X\right)\longrightarrow
		\ttop\left(U,\P_{\alpha,\gamma}X\right)\]
		induced by the
		composition law of $X$.
	\end{enumerate}
\end{nota}

The flow $\{U,X\}_S$ is functorial with respect to $U$ and $X$ (contravariant with respect to $U$ and covariant with respect to $X$). The flow $\{\varnothing,Y\}_S$ is the flow having the same set of states as $Y$ and exactly one non-constant execution path between two points of $Y^0$. The flow $\{\{0\},X\}_S$ is canonically isomorphic to $X$ for all flows $X$. Maybe the latter assertion deserves a little explanation because it is precisely why cartesian closedness matters (the fact that $\{0\}$ is exponentiable is actually sufficient). For all topological spaces $U$, we have the natural bijections \[\top(U,\P_{\alpha,\beta}X) \iso\top(U\p \{0\},\P_{\alpha,\beta}X) \iso \top(U,\ttop(\{0\},\P_{\alpha,\beta}X)).\] Thus by Yoneda, we obtain the homeomorphism \[\P_{\alpha,\beta}X\iso\ttop(\{0\},\P_{\alpha,\beta}X)).\]

\begin{nota}  Let $n\geq 1$. Denote by $\mathbf{D}^n = \{b\in \mathbb{R}^n, |b| \leq 1\}$ the $n$-dimensional disk, and by $\mathbf{S}^{n-1}$ the $(n-1)$-dimensional sphere. By convention, let $\mathbf{D}^{0}=\{0\}$ and $\mathbf{S}^{-1}=\varnothing$. Let $I^{gl}_+=\{\glob(\mathbf{S}^{n-1}) \subset \glob(\mathbf{D}^{n})\mid n\geq 0\} \cup \{C:\varnothing \to \{0\},R:\{0,1\} \to \{0\}\}$. \end{nota}

We recall an elementary lemma about the model category $\top$ which is a straightforward consequence of the fact that the Quillen model structure is Cartesian monoidal: 

\begin{lem} \label{fib}
	Let $i:U\to V$ be a cofibration of $\top$. Then for all topological spaces $X$, the map $i^*:\ttop(V,X)\to \ttop(U,X)$ is a fibration. If moreover $i$ is a weak equivalence, then the map $i^*:\ttop(V,X)\to \ttop(U,X)$ is a trivial fibration. 
\end{lem}

We can now easily carry out the construction of the model category structure.

\bth \label{simplification} There exists a unique model category structure such that $I^{gl}_+$ is the set of generating cofibrations and such that all objects are fibrant. \eth

\bpf We have to check the hypotheses of Theorem~\ref{Isa}. The category $\dtop$ is locally presentable since $\top$ is locally presentable (see for example the proof of \cite[Proposition~6.11]{nonexistence}). That all $\glob(\mathbf{S}^{n-1})$ for all $n\geq 0$ are $I^{gl}_+$-cofibrant comes from the fact that the $(n-1)$-sphere is cofibrant in $\top$. We can factor the relative codiagonal map $\mathbf{D}^{n} \sqcup_{\mathbf{S}^{n-1}}\mathbf{D}^{n} \to \mathbf{D}^{n}$ as a composite $\mathbf{D}^{n} \sqcup_{\mathbf{S}^{n-1}}\mathbf{D}^{n} \subset \mathbf{D}^{n+1} \to \mathbf{D}^{n}$ for all $n\geq 0$. Thus for $U\to V$ being one of the maps $\glob(\mathbf{S}^{n-1}) \subset \glob(\mathbf{D}^{n})$ for $n\geq 0$, we set $C_U(V) = \glob(\mathbf{D}^{n+1})$. We have a pushout diagram of topological spaces 
\[\xymatrix
{
	\mathbf{S}^{n}\fr{}\fd{} & \mathbf{D}^{n} \sqcup_{\mathbf{S}^{n-1}}\mathbf{D}^{n}\fd{}\\
	\mathbf{D}^{n+1}\fr{} & \cocartesien\mathbf{D}^{n+1}
}\]
which gives rise to the pushout diagram of flows 
\[\xymatrix
{
	\glob(\mathbf{S}^{n})\fd{}\fr{} & \glob(\mathbf{D}^{n}) \sqcup_{\glob(\mathbf{S}^{n-1})}\glob(\mathbf{D}^{n})\fd{}\\
	\glob(\mathbf{D}^{n+1})\fr{} & \cocartesien \glob(\mathbf{D}^{n+1})
}\] 
for all $n\geq 0$ by Proposition~\ref{connected-colim-glob}. This implies that for $U\to V$ being one of the maps $\glob(\mathbf{S}^{n-1}) \subset \glob(\mathbf{D}^{n})$ for $n\geq 0$, the map $V\sqcup_U V\to C_U(V)$ belongs to $\cell(I^{gl}_+)$. The map $C:\varnothing\to\{0\}$ gives rise to the relative codiagonal map $\{0\}\sqcup \{0\}\to \{0\}$. Thus we set $C_\varnothing(\{0\})=\{0\}$. In this case, the map $V\sqcup_U V\to C_U(V)$ is $R:\{0,1\}\to\{0\}$ which belongs to $\cell(I^{gl}_+)$. The map $R:\{0,1\}\to \{0\}$ gives rise to the relative codiagonal map $\id_{\{0\}}$. Thus we set $C_{\{0,1\}}(\{0\})=\{0\}$. In this case, the map $V\sqcup_U V\to C_U(V)$ is $\id_{\{0\}}$ which belongs to $\cell(I^{gl}_+)$. The set of generating trivial cofibrations will be therefore the set of maps  $\glob(\mathbf{D}^{n}) \subset \glob(\mathbf{D}^{n+1})$ for $n\geq 0$. Let $\cocyl(X)=\{[0,1],X\}_S$ for all flows $X$. The composite map $\{0,1\}\subset [0,1] \to \{0\}$ yields a natural composite map of flows $X\iso \{\{0\},X\}_S \to \cocyl(X) \to \{\{0,1\},X\}_S$ which is constant on states and which gives rise to the composite continuous map $\P_{\alpha,\beta}X \to \ttop([0,1],\P_{\alpha,\beta}X) \to \P_{\alpha,\beta}X \p\P_{\alpha,\beta}X$ on the spaces of paths for all $(\alpha,\beta)\in X^0\p X^0$. We obtain a natural composite map of flows $X\stackrel{\tau}\longrightarrow \cocyl(X) \stackrel{\pi}\longrightarrow X\p X$ since the set of states of $X\p X$ is $X^0 \p X^0$ and the space of paths from $(\alpha,\alpha')$ to $(\beta,\beta')$ is $\P_{\alpha,\beta}X \p \P_{\alpha',\beta'}X$ by \cite[Theorem~4.17]{model3}. We have obtained a path object in the sense of Theorem~\ref{Isa}. Since the maps $\pi_0$ and $\pi_1$ are bijective on states, they satisfy the RLP with respect to $\{C:\varnothing \to \{0\},R:\{0,1\} \to \{0\}\}$ by Proposition~\ref{ortho2}. By Proposition~\ref{ortho3}, the maps $\pi_0$ and $\pi_1$ satisfy the RLP with respect to $\glob(\mathbf{S}^{n-1}) \subset \glob(\mathbf{D}^{n})$ for $n\geq 0$ if and only if the evaluation maps $\ttop([0,1],\P_{\alpha,\beta}X) \rightrightarrows \P_{\alpha,\beta}X$  on $0$ and $1$ satisfy the RLP with respect to the inclusion $\mathbf{S}^{n-1}\subset \mathbf{D}^{n}$ for $n\geq 0$ and for all $(\alpha,\beta)\in X^0\p X^0$, i.e. if and only if the evaluation maps $\ttop([0,1],\P_{\alpha,\beta}X) \rightrightarrows \P_{\alpha,\beta}X$ are trivial fibrations for all $(\alpha,\beta)\in X^0\p X^0$. The latter fact is a consequence of Lemma~\ref{fib} and from the fact that the inclusions $\{0\}\subset [0,1]$ and $\{1\}\subset [0,1]$ are trivial cofibrations of $\top$. Finally we have to check that the map $\pi:\cocyl(X)\to X\p X$ satisfies the RLP with respect to the maps $\glob(\mathbf{D}^{n}) \subset \glob(\mathbf{D}^{n+1})$ for $n\geq 0$. By Proposition~\ref{ortho3} again, it suffices to prove that the map $\ttop([0,1],\P_{\alpha,\beta} X) \to \ttop(\{0,1\},\P_{\alpha,\beta} X) = \P_{\alpha,\beta}X\p \P_{\alpha,\beta}X$ is a fibration of topological spaces for all $(\alpha,\beta)\in X^0\p X^0$. By Lemma~\ref{fib} again, this comes from the fact that the inclusion $\{0,1\}\subset [0,1]$ is a cofibration of $\top$.
\epf

This model category structure coincides with the one of \cite{model3}.

\section{Left determinedness of the model category of flows}
\label{section3}

Let us now recall the definition of a left determined model category: 

\bd \label{localizer}
Let $I$ be a set of maps of a locally presentable category $\K$. A class of maps $\W$ is a {\rm localizer (with respect to $I$)} or an {\rm $I$-localizer} if $\W$ satisfies: 
\begin{itemize}
	\item Every map satisfying the RLP with respect to the maps of $I$ belongs to $\W$.
	\item $\W$ is closed under retract and satisfies the $2$-out-of-$3$ property.
	\item The class of maps $\cof(I)\cap \W$ is closed under pushout and transfinite composition.
\end{itemize}
The class of all maps is an $I$-localizer. The class of $I$-localizers is closed under arbitrarily large intersection. Therefore there exists a smallest $I$-localizer denoted by $\W_I$.
\ed

\bd \cite{rotho} A combinatorial model category $\K$ with the set of generating cofibrations $I$ is {\rm left determined} if the class of weak equivalences is $\W_I$. \ed

Consider a combinatorial model category $\K$ such that all objects are fibrant with a class of weak equivalences $\W$ and a set of generating cofibrations $I$. The localizer $\W_I$ could be strictly smaller than $\W$. If $\W_I$ is the class of weak equivalences of a model category structure on $\K$, then all objects of this model category structure are fibrant, and therefore $\W=\W_I$. To the best of our knowledge, we can only say, using \cite[Theorem~2.2]{rotho}, that every combinatorial model category such that all objects are fibrant is left determined if we assume Vop\v{e}nka's principle. Since in any model category, two fibrant objects are weakly equivalent if and only if they are related by a span of trivial fibrations, and since all trivial fibrations belong to the smallest localizer, it is also true that there is an equivalence of categories $\K[\W_I^{-1}]\simeq \K[\W^{-1}]$ between the categorical localizations of $\K$ with respect to $\W_I$ and $\W$ if all objects of the combinatorial model category $\K$ are fibrant. Note that in \cite{Isaev}, a localizer is just the class of weak equivalences of a model category structure. In the latter sense, a model category of fibrant objects has a minimal localizer indeed. 

In our case, it is possible to conclude that the model category is left determined without assuming Vop\v{e}nka's principle by adapting a technique we learned from Marc Olschok for the model category of topological spaces.

\bth \label{leftdetflow} The model category of flows is left determined. \eth

\bpf Let $f:X\to Y$ be a weak equivalence of flows. Then $f$ factors as a composite $f=f_2.f_1$ where $f_1$ is a trivial cofibration, i.e.  $f_1\in \cof(\{\glob(\mathbf{D}^{n}) \subset \glob(\mathbf{D}^{n+1})\mid n\geq 0\})$ and where $f_2$ is a trivial fibration. In particular $f_2$ satisfies the RLP with respect to $C:\varnothing\to \{0\}$ and $R:\{0,1\}\to \{0\}$. Thus $f_2$ is bijective on states by Proposition~\ref{ortho2}. The functor $X\mapsto X^0$ from $\dtop$ to $\set$ is colimit-preserving since it has a right adjoint (the functor taking a set $S$ to the flow with the set of states $S$ and exactly one path between each pair of states). Therefore $f_1$ is bijective on states since the maps $\glob(\mathbf{D}^{n}) \subset \glob(\mathbf{D}^{n+1})$ for all $n\geq 0$ are bijective on states. We deduce that $f$ is bijective on states. Consider the commutative diagram of flows
\[
\xymatrix
{
	X \fR{\exists !\overline{f}} \fD{f} && N_f \cartesien \fR{p'} \fD{f'} && X \fD{f} \\
	&& && \\
	Y \fR{\tau} && \{[0,1],Y\}_S \fR{\pi_0} \fD{\pi_1} && Y \\
	&& && \\
	&& Y &&
}
\]
where the existence of $\overline{f}$ comes from the universal property of the pullback. All arrows are bijective on states. Using the fact that the functor $\P:\dtop\to\top$ is limit-preserving by \cite[Theorem~4.17]{model3}, one obtains the commutative diagram  of topological spaces 
\[
\xymatrix
{
	\P X \fR{\P \overline{f}} \fD{\P f} && \P N_f \cartesien \fR{\P p'} \fD{\P f'} && \P X \fD{\P f} \\
	&& && \\
	\P Y \fR{\P\tau} && \ttop([0,1],\P Y) \fR{\P(\pi_0)} \fD{\P(\pi_1)} && \P Y \\
	&& && \\
	&& \P Y &&
}
\]
By \cite[Proposition~4.64]{MR1867354}, the map $\P(\pi_1).\P(f')$ is a Hurewicz fibration, and therefore a fibration of the model category of $\top$. By Proposition~\ref{ortho3}, $\pi_1.f'$ satisfies the RLP with respect to all trivial cofibrations of flows, i.e. $\pi_1.f'$ is a fibration of flows. The maps $\P(\pi_0)$ and $\P(\pi_1)$ are trivial fibrations of the model category of $\top$ by Lemma~\ref{fib}.  By Proposition~\ref{ortho3} and Proposition~\ref{ortho2}, $\pi_0$ and $\pi_1$ satisfy the RLP with respect to all cofibrations of $\dtop$, i.e. $\pi_0$ and $\pi_1$ are trivial fibrations of flows. Thus $p'$ is a trivial fibration of $\dtop$ since it is a pullback of a trivial fibration. Since $f$ is a weak equivalence of flows by hypothesis, we deduce by the $2$-out-of-$3$ property that $f'$ is a weak equivalence of $\dtop$. Thus $\pi_1.f'$ is a trivial fibration of flows as well. We have $p'.\overline{f}=\id_X$. Since $p'$ is a trivial fibration, it belongs to the smallest localizer. Therefore by the $2$-out-of-$3$ property, $\overline{f}$ belongs to the smallest localizer. Since $\pi_1.f'$ is a trivial fibration, it belongs to the smallest localizer as well. Since $\pi_1.f'.\overline{f} = \pi_1.\tau.f=f$, we deduce that $f$ belongs to the smallest localizer.
\epf

\section{Concluding remarks}
\label{conc}

The hypothesis that $\top$ is locally presentable can be removed. Theorem~\ref{simplification} and Theorem~\ref{leftdetflow} hold by working in any bicomplete cartesian closed full subcategory of the general category of topological spaces containing all CW-complexes. But then, we have to check that all domains and all codomains of the maps of $I_{gl}^+$ are small relative to $\cell(I_{gl}^+)$. This is done in \cite[Section~11]{model3} and there is no known way to avoid the use of some difficult topological arguments. However, the model category of flows is left proper but not cellular because of the presence of $R:\{0,1\}\to \{0\}$ in the generating cofibrations. So outside the framework of locally presentable categories, we have no tools to prove the existence of any homotopical localization and to study the homotopical localization of $\dtop$ with respect to the refinement of observation.

\end{document}